\documentclass[a4paper,twoside,11pt]{article} 
 
\usepackage{amsmath, amsthm, amssymb, graphicx, txfonts, bm}

\usepackage{setspace}
\usepackage[inner=2.4cm,outer=2.4cm,top=2.4cm,bottom=2.4cm,marginparwidth=90pt]{geometry} 

\usepackage{accents}
\usepackage{times}

\usepackage{tikz} 
\usepackage{pgf}
\usetikzlibrary{positioning}

\numberwithin{equation}{section} 

\usepackage{graphics}
\usepackage{mathrsfs}
\usepackage{framed}
\usepackage{comment}

\newcommand{\Rone}{\mathbb{R}}

\newcommand{\Comone}{\mathbb{C}}
\newcommand{\SR}{\mathscr{S}_R^n}
\newcommand{\Rnone}{\mathbb{R}^{n+1}}

\renewcommand{\proof}{\smallskip\par {\noindent\textbf{Proof:}}\quad}

\newcommand{\proofoftheoremone}{\smallskip\par {\noindent\textbf{Proof of Theorem 1.1 :\quad}}}

\newtheorem{theorem}{Theorem}[section]
\newtheorem{lemma}[theorem]{Lemma}

\newtheorem{corollary}[theorem]{Corollary}

\begin{document}

\title{Motion by Mixed Volume Preserving Curvature Functions Near Spheres}

\author{David Hartley}

\date{\today}

\onehalfspacing

\maketitle

\begin{abstract}
In this paper we investigate the flow of surfaces by a class of symmetric functions of the principal curvatures with a mixed volume constraint. We consider compact surfaces without boundary that can be written as a graph over a sphere. The linearisation of the resulting fully nonlinear PDE is used to prove a short time existence theorem for a large class of surfaces that are sufficiently close to a sphere and, using center manifold analysis, the stability of the sphere as a stationary solution to the flow is determined. We will find that for initial surfaces sufficiently close to a sphere, the flow will exist for all time and converge exponentially to a sphere. This result was shown for the case where the symmetric function is the mean curvature and the constraint is on the $(n+1)$-dimensional enclosed volume by Escher and Simonett \cite{Escher98A}.
\end{abstract}

\section{Introduction}
Given a sufficiently smooth, compact without boundary initial hypersurface $\Omega_0=\bm{X}_0(M^n)\subset\Rnone$ we are interested in finding a family of embeddings $\bm{X}:M^n\times[0,T)\rightarrow\Rnone$ such that
\begin{equation}
\frac{\partial\bm{X}}{\partial t}=\left(h_k-F(\bm{\kappa})\right)\bm{\nu}_{\Omega_t},\ \ \ \bm{X}(\cdot,0)=\bm{X}_0,\ \ \ h_k=\frac{1}{\int_{M^n}E_{k+1}\,d\mu_t}\int_{M^n} F(\bm{\kappa})E_{k+1}\,d\mu_t,
\label{FullFlow}
\end{equation}
where $\bm{\kappa}=\left(\kappa_1,\ldots,\kappa_n\right)$ and $\kappa_i$ are the principal curvatures of the surface $\Omega_t=\bm{X}(M^n,t)=\bm{X}_t(M^n)$, $\bm{\nu}_{\Omega_t}$ and $\,d\mu_t$ are the outward pointing normal and induced measure of $\Omega_t$ respectively and $k$ is a fixed integer between $-1$ and $n-1$. Here $E_l$ denotes the $l^{th}$ elementary symmetric function of the principal curvatures
\begin{equation*}
E_{l}=\left\{\begin{array}{ll} 1 & l=0\\ \sum_{1\leq i_1<\ldots<i_l\leq n} \kappa_{i_1}\kappa_{i_2}\ldots\kappa_{i_l} & l=1,\ldots,n.\end{array}\right.
\end{equation*}
The function $F(\bm{\kappa})$ is a smooth symmetric function with $\frac{\partial F}{\partial \kappa_i}(\bm{\kappa}_0)>0$, where $\bm{\kappa}_0=\left(\frac{1}{R},\ldots,\frac{1}{R}\right)$, for some $R\in\Rone^+$. For a fixed $k$ the flow can be seen to preserve a certain quantity, which for convex hypersurfaces is the $(n-k)^{th}$ mixed volume (see Section \ref{Sec2}).

This flow has been studied previously in \cite{McCoy05} where it was proved that under some additional conditions on $F$, for example homogeneity of degree one and convexity or concavity, initially convex hypersurfaces admit a solution for all time and that the hypersurfaces converge to a sphere as $t\to\infty$. This result had previously been proved for the specific case where $F(\bm{\kappa})=H$, the mean curvature, in \cite{McCoy04} and, if in addition, $k=-1$ (in which the flow is the well known volume preserving mean curvature flow) in \cite{Huisken87}. Other results for the volume preserving mean curvature flow include average mean convex hypersurfaces with initially small traceless second fundamental form converging to spheres (see \cite{Li09}) and hypersurfaces that are graphs over spheres with a height function close to zero, in a certain function space, converging to spheres (see \cite{Escher98A}). Techniques similar to those in this paper were used to study volume preserving mean curvature flow for hypersurfaces close to a cylinder \cite{Hartley12}.

When $k=-1$ and $F(\bm{\kappa})=H_m^{\beta}$, with $m\beta>1$ and $H_m=\binom{n}{m}^{-1}E_m$ is the $m^{th}$ mean curvature, the flow has been shown to take initially convex hypersurfaces that satisfy a pinching condition to spheres; the pinching condition is of the form $E_n>CH^n>0$, where $C$ is a constant depending on the parameters of the flow \cite{Cabezas10}. When $m=1$ the flow is the volume preserving analogue of the powers of mean curvature flow introduced in \cite{Schulze05}.

The main result of this paper is
\begin{theorem}
If $\Omega_0$ is a graph over the sphere $\SR$ with height function sufficiently small in $C^{2+\alpha}\left(\SR\right)$, $0<\alpha<1$, then its flow by (\ref{FullFlow}) exists for all time and converges exponentially fast to a sphere as $t\rightarrow\infty$, with respect to the $C^{2+\alpha}\left(\SR\right)$ topology.
\label{MainR}
\end{theorem}

This result is the analogue of part c) of the Main Result in \cite{Escher98A} for fully nonlinear equations. Some differences include the fact that the smoothness of the hypersurface after the initial time is not guaranteed and hence the convergence is only with respect to the $C^{2+\alpha}\left(\SR\right)$ topology, instead of the convergence with respect to the $C^l$ topology proved in \cite{Escher98A}. This is because volume preserving mean curvature flow is quasilinear while the flows in this paper are in general fully nonlinear. Also control of the curvature is required so the initial hypersurface is small in $C^{2+\alpha}\left(\SR\right)$.

In Section \ref{Sec2} of this paper we convert the flow (\ref{FullFlow}) to a PDE for the graph function and also introduce the spaces and notation that will be used throughout the paper. The section ends with a lemma that the flow preserves a certain mixed volume. In Section \ref{Sec3} we consider the problem as an ODE on Banach spaces and determine the linearisation of the speed. This leads to a short-time existence theorem for the flow that includes some initially non-smooth hypersurfaces. In the final section the eigenvalues of the linearised operator are determined and a center manifold is constructed. The proof of the main result is finished by showing that the center manifold consists entirely of spheres and is exponentially attractive.

The $(n-k)^{th}$-mixed volumes, for $k\geq1$, are only well defined for convex hypersurfaces (see \cite{Andrews01}). In this paper we do not make the explicit assumption that the hypersurfaces are convex, however with the closeness, in $C^{2+\alpha}\left(\SR\right)$, to a sphere condition this may well be the case. We will continue to refer to the flow as mixed volume preserving with the understanding that it preserves a specific quantity that coincide with a mixed volume for convex hypersurfaces.

The author would like to thank Maria Athanassenas for her support, advice and help in preparing this paper, Todd Oliynyk for his advice and encouragement in dealing with center manifolds and to Monash University and the School of Mathematical Sciences for their support.

\section{Notation and Preliminaries}
\label{Sec2}
In this paper we consider $M^n=\SR$, a given sphere of radius $R$, and $\bm{X}_0(\bm{p})=\bm{p}+\rho_0(\bm{p})\bm{\nu}_{\SR}(\bm{p})$, $\bm{p}\in\SR$, so that $\bm{X}_0$ is a graph over $\SR$. The volume form on such a hypersurface will be denoted by $\,d\mu_{\rho}$ and we let $\mu_{\rho}$ be the function such that $\,d\mu_{\rho}=\mu_{\rho}\,d\mu_0$. We now proceed as Escher and Simonett \cite{Escher98A} and convert the flow to an evolution equation for the height function $\rho:\SR\times[0,T)\rightarrow\Rone$. Up to a tangential diffeomorphism the flow (\ref{FullFlow}) is equivalent to
\begin{equation}
\frac{\partial\rho}{\partial t}=\sqrt{1+\frac{R^2}{(R+\rho)^2}\left|\nabla\rho\right|^2}\left(h_{k,\rho}-F(\bm{\kappa}_{\rho})\right),\ \ \ \ \rho(\cdot,0)=\rho_0,
\label{GraphFlow}
\end{equation}
where $\bm{\kappa}_{\rho}$ is the principal curvature vector of the hypersurface defined by $\rho(\cdot,t)$ and $\nabla$ denotes the gradient on $\SR$ (see \cite{Bode07}).

The graph functions $\rho$ are chosen in the H\"older spaces, $C^{l,\alpha}\left(\SR\right)$, for $\alpha\in(0,1)$, $l\in\mathbb{N}$.These spaces are the interpolation spaces between the $C^l$ spaces (see \cite{Lunardi95}),
\begin{equation*}
C^{l\theta}\left(\SR\right)=\left(C\left(\SR\right),C^l\left(\SR\right)\right)_{\theta,\infty},
\end{equation*}
where $(\cdot,\cdot)_{\theta,\infty}$ is an interpolation functor for each $\theta\in(0,1)$ and is defined for $Y\subset X$ as follows:
\begin{equation*}
(X,Y)_{\theta,\infty}=\left\{x\in X:\lim_{t\to0}t^{-\theta}K(t,x,X,Y)<\infty\right\},\ \ \ K(t,x,X,Y)=\inf_{a\in Y}\left(\left\|x-a\right\|_{X}+t\left\|a\right\|_Y\right).
\end{equation*}
By Corollary 1.2.18 in \cite{Lunardi95} we have
\begin{equation}
\left(C^{\alpha}\left(\SR\right),C^{2+\alpha}\left(\SR\right)\right)_{\theta,\infty}= C^{\alpha+2\theta}\left(\SR\right),\ \ \ \alpha+2\theta\notin\mathbb{N},
\label{Interp}
\end{equation}
where $\alpha\in(0,1)$.

For an operator between function spaces $G:Y\rightarrow \tilde{Y}$ we denote the Fr\'echet derivative by $\partial G$. A linear operator, $A:Y\subset X\rightarrow X$, is called \textit{sectorial} if there exist $\theta\in\left(\frac{\pi}{2},\pi\right)$ and $\omega\in\Rone$ and $M>0$ such that
\begin{equation*}
\left\{\begin{array}{ll} (i) & \rho(A)\supset S_{\theta,\omega}=\{\lambda\in\Comone:\lambda\neq\omega,||\arg(\lambda-\omega)|<\theta\},\\ (ii) & \|R(\lambda,A)\|_{L\left(X\right)}\leq\frac{M}{|\lambda-\omega|}\textrm{ for all }\lambda\in S_{\theta,\omega}.\end{array}\right.
\end{equation*}
Here $\rho(A)$  is the resolvent set, $R(\lambda,A)=(\lambda I-A)^{-1}$ is the resolvent operator and $\|\cdot\|_{L(X)}$ is the standard linear operator norm (see \cite{Lunardi95}).

For all hypersurfaces $\Omega=X(M)$ we define the quantity
\begin{equation*}
V_{n-k}\left(\Omega\right) = \left\{\begin{array}{ll} \mathrm{Vol}(\Phi) & k=-1\\ \left((n+1)\binom{n}{k}\right)^{-1}\int_{M}E_k\,d\mu & k=0,\ldots,n-1,\end{array}\right.
\end{equation*}
where $\Phi$ is the $(n+1)$-dimensional region contained inside $\Omega$, for convex hypersurfaces this agrees with the mixed volumes.

\begin{lemma}
For an initially smooth, compact, convex hypersurface without boundary, $\Omega_0$, the flow (\ref{FullFlow}) preserves the value of $V_{n-k}$, i.e. $V_{n-k}\left(\Omega_t\right)=V_{n-k}\left(\Omega_0\right)$ as long as the flow exists.
\end{lemma}

\proof
This is proved through a calculation of the evolution equations, see \cite{McCoy05}.
$\Box$

\section{Graphs over Spheres}
\label{Sec3}
The flow in equation (\ref{GraphFlow}) can be considered as an ordinary differential equation between Banach spaces. Set $0<\alpha<1$ and define
\begin{equation*}
G:C^{2+\alpha}\left(\SR\right)\rightarrow C^{\alpha}\left(\SR\right),\ \ \ G(\rho):=L_{\rho}\left(h_{k,\rho}-F(\bm{\kappa}_\rho)\right),\ \ \ L_{\rho}:= \sqrt{1+\frac{R^2}{(R+\rho)^2}|\nabla \rho|^2}.
\end{equation*}
The flow (\ref{GraphFlow}) is then rewritten as
\begin{equation}
\rho'(t)=G(\rho(t)),\ \ \ \ \rho(0)=\rho_0\in C^{2+\alpha}\left(\SR\right).
\label{GraphODE}
\end{equation}

\begin{lemma} For the linearisation, $\partial G(0)$, of $G$ it holds
\begin{equation*}
\partial G(0)u = \frac{\partial F}{\partial\kappa_1}(\bm{\kappa}_0)\left(\left(\frac{n}{R^2}+\Delta_{\SR}\right)u-\frac{n}{R^2}\fint_{\SR}u\,d\mu_0\right),
\end{equation*}
for $u\in C^{2+\alpha}\left(\SR\right)$.
\label{LinOp2}
\end{lemma}

Note that only the derivative of $F(\bm{\kappa})$ with respect to $\kappa_1$ appears in this formula for convenience, since $ \frac{\partial F}{\partial\kappa_1}(\bm{\kappa}_0)= \frac{\partial F}{\partial\kappa_i}(\bm{\kappa}_0)$ for all $i=1,\ldots,n$.

\proof
Firstly note that $L_0=1$ and that $\left.\partial L_{\rho}\right|_{\rho=0}=0$ and by linearising the curvature function we find
\begin{align*}
\left.\partial F(\bm{\kappa}_{\rho})\right|_{\rho=0} &= \left.\sum_{i=1}^n\frac{\partial F}{\partial\kappa_i}(\bm{\kappa}_{\rho})\partial \kappa_{i,\rho}\right|_{\rho=0}\\
&=\frac{\partial F}{\partial\kappa_1}(\bm{\kappa}_0)\left.\sum_{i=1}^n\partial \kappa_{i,\rho}\right|_{\rho=0}\\
&=\frac{\partial F}{\partial\kappa_1}(\bm{\kappa}_0)\left.\partial H_{\rho}\right|_{\rho=0}.
\end{align*}
It follows that for $u\in C^{2+\alpha}\left(\SR\right)$
\begin{align*}
\left.\partial h_{k,\rho}\right|_{\rho=0}u &= \left.\partial\left(\frac{1}{\int_{\SR}E_{k+1,\rho}\mu_{\rho}\,d\mu_0}\int_{\SR}E_{k+1,\rho}F(\bm{\kappa}_{\rho})\mu_{\rho}\,d\mu_0\right)\right|_{\rho=0}u\\
&=\frac{1}{\left(\int_{\SR}E_{k+1,0}\,d\mu_0\right)^2}\left(\int_{\SR}E_{k+1,0}\,d\mu_0\left.\partial\left(\int_{\SR}E_{k+1,\rho}F(\bm{\kappa}_{\rho})\mu_{\rho}\,d\mu_0\right)\right|_{\rho=0}u\right.\\
&\hspace{4cm}\left. -\int_{\SR}E_{k+1,0}F(\bm{\kappa}_{0})\,d\mu_0\left.\partial\left(\int_{\SR}E_{k+1,\rho}\mu_{\rho}\,d\mu_0\right)\right|_{\rho=0}u\right)\\
&=\frac{1}{\int_{\SR}E_{k+1,0}\,d\mu_0}\left(\int_{\SR}E_{k+1,0}\left.\partial F(\bm{\kappa}_{\rho})\right|_{\rho=0}u+F(\bm{\kappa}_{0})\left.\partial\left( E_{k+1,\rho}\mu_{\rho}\right)\right|_{\rho=0}u\,d\mu_0\right.\\
&\hspace{4cm}\left. -F(\bm{\kappa}_{0})\int_{\SR}\left.\partial\left( E_{k+1,\rho}\mu_{\rho}\right)\right|_{\rho=0}u\,d\mu_0\right)
\end{align*}
\begin{align*}
\left.\partial h_{k,\rho}\right|_{\rho=0}u&=\frac{\frac{\partial F}{\partial\kappa_1}(\bm{\kappa}_0)}{\left|\SR\right|}\int_{\SR}\left.\partial H_{\rho}\right|_{\rho=0}u\,d\mu_0\\
&=\frac{\partial F}{\partial\kappa_1}(\bm{\kappa}_0)\fint_{\SR}\left.\partial H_{\rho}\right|_{\rho=0}u\,d\mu_0.
\end{align*}
It was shown in \cite{Escher98B} that
\begin{equation*}
\left.\partial H_{\rho}\right|_{\rho=0}=-\left(\frac{n}{R^2}+\Delta_{\SR}\right),
\end{equation*}
so combining these results gives, for $u\in C^{2+\alpha}\left(\SR\right)$,
\begin{equation}
\partial G(0)u = \frac{\partial F}{\partial\kappa_1}(\bm{\kappa}_0)\left(\left(\frac{n}{R^2}+\Delta_{\SR}\right)u-\fint_{\SR}\left(\frac{n}{R^2}+\Delta_{\SR}\right)u\,d\mu_0\right).
\label{LinOp1}
\end{equation}
The divergence theorem gives the result.
$\Box$

\begin{lemma}
For any $\alpha_0$ such that $0<\alpha_0<\alpha$ there exists a neighbourhood, $O_1$, of $0\in C^{2+\alpha}\left(\SR\right)$ such that the operator $\partial G(\rho)$ is the part in $C^{\alpha}\left(\SR\right)$ of a sectorial operator $A_{\rho}:C^{2+\alpha_0}\rightarrow C^{\alpha_0}\left(\SR\right)$ for all $\rho\in O_1$.
\label{SecLem}
\end{lemma}

\proof
We set $\bar{G}:C^{2+\alpha_0}\left(\SR\right)\rightarrow C^{\alpha_0}\left(\SR\right)$ with $\bar{G}(\rho):=L_{\rho}\left(h_{k,\rho}-F\left(\bm{\kappa}_{\rho}\right)\right)$ so that with $A_{\rho}=\partial \bar{G}(\rho)$ it is clear that $\partial G(\rho)$ is the part in $C^{\alpha}\left(\SR\right)$ of $A_{\rho}$. It remains to show that there exists $O_1$ such that $A_{\rho}$ is sectorial for $\rho\in O_1$.

As $\left.\partial H_{\rho}\right|_{\rho=0}=-\left(\frac{n}{R^2}+\Delta_{\SR}\right)$ is a uniformly elliptic operator, its negative is sectorial. Now the operator $A_0:C^{2+\alpha_0}\left(\SR\right)\rightarrow C^{\alpha_0}\left(\SR\right)$, defined by
\begin{equation*}
A_0u=\left(\frac{n}{R^2}+\Delta_{\SR}\right)u-\frac{n}{R^2}\fint_{\SR}u\,d\mu_0,
\end{equation*}
can be seen to be sectorial: using the definition of sectorial, there exists $M>0$, $\theta\in\left(\frac{\pi}{2},\pi\right)$ and $\omega\in\Rone$ such that
\begin{equation*}
\left|\lambda-\omega\right|\left\|u\right\|_{C^{2+\alpha_0}\left(\SR\right)}\leq M\left\|\left(\lambda+\left.\partial H_{\rho}\right|_{\rho=0}\right)u\right\|_{C^{\alpha_0}\left(\SR\right)},\ \ \textrm{for all } \lambda\in S_{\theta,\omega}\textrm{ and } u\in C^{2+\alpha_0}\left(\SR\right).
\end{equation*}
Therefore
\begin{align*}
\left\|(\lambda-A_0)u\right\|_{C^{\alpha_0}\left(\SR\right)}&=\left\|\left(\lambda-\left(\frac{n}{R^2}+\Delta_{\SR}\right)\right)u-\left(-\frac{n}{R^2}\fint_{\SR}u\,d\mu_0\right)\right\|_{C^{\alpha_0}\left(\SR\right)}\\
&\geq\left\|\left(\lambda+\left.\partial H_{\rho}\right|_{\rho=0}\right)u\right\|_{C^{\alpha_0}\left(\SR\right)}-\left\|-\frac{n}{R^2}\fint_{\SR}u\,d\mu_0\right\|_{C^{\alpha_0}\left(\SR\right)}\\
&\geq\frac{|\lambda-\omega|}{M}\|u\|_{C^{2+\alpha_0}\left(\SR\right)}-\frac{n}{R^2}\|u\|_{C(\SR)}\\
&\geq\frac{|\lambda-\omega|}{M}\|u\|_{C^{2+\alpha_0}\left(\SR\right)}-\frac{n}{R^2}\|u\|_{C^{2+\alpha_0}\left(\SR\right)}\\
&=|\lambda|\left(\frac{\left|1-\frac{\omega}{\lambda}\right|}{M}-\frac{n}{R^2|\lambda|}\right)\|u\|_{C^{2+\alpha_0}\left(\SR\right)}
\end{align*}
By the reverse triangle inequality and by taking $Re(\lambda)\geq\omega_1:=\left(2|\omega|+\frac{2Mn}{R^2}\right)$ we get
\begin{align*}
\left\|(\lambda-A_0)u\right\|_{C^{\alpha_0}\left(\SR\right)}&\geq|\lambda|\left(\frac{1}{M}-\frac{|\omega|}{M|\lambda|}-\frac{n}{R^2|\lambda|}\right)\|u\|_{C^{2+\alpha_0}\left(\SR\right)}\\
&\geq\frac{|\lambda|}{2M}\|u\|_{C^{2+\alpha_0}\left(\SR\right)},
\end{align*}
Thus, by applying Proposition 2.1.11 in \cite{Lunardi95}, $A_0$ is sectorial. This then implies by Proposition 2.4.2 in \cite{Lunardi95} that $A_{\rho}= \partial\bar{G}(0)+\left(\partial\bar{G}(\rho)-\partial\bar{G}(0)\right)$ is sectorial for all $\rho$ in a neighbourhood of zero, $O_2\subset C^{2+\alpha_0}\left(\SR\right)$, so that $O_1=O_2\cap C^{2+\alpha}\left(\SR\right)$.
$\Box$

\begin{theorem}
There are constants $\delta,r>0$ such that if $\left\|\rho_0\right\|_{C^{2+\alpha}\left(\SR\right)}\leq r$ then equation (\ref{GraphODE}) has a solution
\begin{equation*}
\rho\in C\left([0,\delta],C^{2+\alpha}\left(\SR\right)\right)\cap C^1\left([0,\delta],C^{\alpha}\left(\SR\right)\right)\textrm{ with }\rho(0)=\rho_0.
\end{equation*}
\end{theorem}

\proof
This existence theorem is a result of Theorem 8.4.1 in \cite{Lunardi95} by setting $\bar{t}=t_0=0$ and $\bar{u}=0$. In order to satisfy the assumption of the theorem it must be shown that there exists a neighbourhood of zero, $O\subset C^{2+\alpha}\left(\SR\right)$, such that $G$ and $\partial G$ are continuous on $O$ and that for every $\bar{\rho}\in O$ the operator $\partial G(\bar{\rho})$ is the part in $C^{\alpha}\left(\SR\right)$ of a sectorial operator $A:C^{2+\alpha_0}\left(\SR\right)\rightarrow C^{\alpha_0}\left(\SR\right)$.

As in \cite{Andrews09} Remark 1, since $F$ is a smooth symmetric function of the principal curvatures it is also a smooth function of the elementary symmetric functions, which depend smoothly on the components of the Weingarten map. It is easily seen that the Weingarten map depends smoothly on $\rho\in C^{2+\alpha}\left(\SR\right)$ so that $G$ depends smoothly on $\rho\in C^{2+\alpha}\left(\SR\right)$ inside a neighbourhood, $O_3$, where if $\rho\in O_3$ we have $\int_{\SR} E_{k+1,\rho}\,d\mu_{\rho}>0$ and $\rho(\bm{p})>-R$ for all $\bm{p}\in\SR$ (note if $k=-1$ the former is always satisfied). The sectorial condition was established in Lemma \ref{SecLem} for a neighbourhood $O_1$, so the proof is complete by setting $O=O_3\cap O_1$.
$\Box$

\section{Stability around Spheres}
As we are considering the flow locally about $\rho=0$, it is convenient to rewrite (\ref{GraphODE}) highlighting the dominant linear part
\begin{equation}
\rho'(t)=\partial G(0)\rho(t)+\tilde{G}\left(\rho(t)\right),\ \ \ \ \tilde{G}\left(u\right):=G\left(u\right)-\partial G(0)u.
\label{GraphLin}
\end{equation}
\begin{lemma}
The spectrum $\sigma\left(\partial G(0)\right)$ of $\partial G(0)$ consists of a sequence of isolated non-positive eigenvalues where the multiplicity of the $0$ eigenvalue is $n+2$.
\label{Spectrum}
\end{lemma}

\proof
This follows from \cite{Escher98A} as $\partial G(0)$ is a positive constant multiple of the linear operator in their paper. To be exact, we calculate all the elements of the spectrum. Since $C^{2+\alpha}\left(\SR\right)$ is compactly embedded in $C^{\alpha}\left(\SR\right)$, the spectrum consists entirely of eigenvalues. To characterise the spectrum we first look at the spectrum of the $L^2$-self adjoint operator:
\begin{equation*}
\tilde{A}u=\frac{\partial F}{\partial\kappa_1}(\bm{\kappa}_0)\left(\frac{n}{R^2}+\Delta_{\SR}\right)u.
\end{equation*}
The eigenvalues of the spherical Laplacian are well known to be $\frac{-l(l+n-1)}{R^2}$ for $l\in\mathbb{N}\cup\{0\}$ with eigenfunctions the spherical harmonics of order $l$, denoted by $Y_{l,p}$, $1\leq p\leq M_l$, where
\begin{equation*}
M_l=\binom{l+n}{n}-\binom{l+n-2}{n}.
\end{equation*}
Therefore the eigenfunctions of $\tilde{A}$ are also the spherical harmonics with eigenvalues
\begin{equation*}
\xi_l=\frac{\partial F}{\partial\kappa_1}(\bm{\kappa}_0)\left(\frac{n}{R^2}-\frac{l(l+n-1)}{R^2}\right)=-\frac{\partial F}{\partial\kappa_1}(\bm{\kappa}_0)\frac{(l-1)(l+n)}{R^2}.
\end{equation*}

Returning to the spectrum of $\partial G(0)$, $Y_{0,1}=1$ is still an eigenfunction but with eigenvalue $\lambda_0=0$. The operator $\partial G(0)$ is self adjoint with respect to the $L^2$ inner product on $C^{2+\alpha}\left(\SR\right)$. To see this, consider $u,w\in C^{2+\alpha}\left(\SR\right)$

\begin{align*}
\int_{\SR}\left(\partial G(0)u\right)w\,d\mu_0 & =\int_{\SR}\left(\tilde{A}u-\frac{n}{R^2}\fint_{\SR}u\,d\mu_0\right)w\,d\mu_0\\
&=\int_{\SR}\left(\tilde{A}u\right)w\,d\mu_0-\frac{n}{R^2}\fint_{\SR}u\,d\mu_0\int_{\SR}w\,d\mu_0\\
&=\int_{\SR}u\left(\tilde{A}w\right)\,d\mu_0-\frac{n}{R^2}\int_{\SR}u\,d\mu_0\fint_{\SR}w\,d\mu_0,
\end{align*}
where we used that $\tilde{A}$ is self adjoint since it is a multiple of the Laplacian on the sphere plus a constant.
\begin{align*}
\int_{\SR}\left(\partial G(0)u\right)w\,d\mu_0&=\int_{\SR}u\left(\tilde{A}w-\frac{n}{R^2}\fint_{\SR}w\,d\mu_0\right)\,d\mu_0\\
&=\int_{\SR}u\left(\partial G(0)w\right)\,d\mu_0.
\end{align*}
Therefore we need only consider other eigenfunctions, orthogonal to $Y_{0,1}=1$, in order to characterise the spectrum. This means that for an eigenfunction $u$
\begin{equation*}
\int_{\SR}u\,d\mu_0=0,
\end{equation*}
hence by Lemma \ref{LinOp2}, $\partial G(0)u=\tilde{A}u$. The remaining eigenfunctions of $\partial G(0)$ are then the remaining eigenfunctions of $\tilde{A}$ , with the same eigenvalues. So the spectrum of $\partial G(0)$ consists of the eigenvalues
\begin{equation*}
\lambda_l=\left\{\begin{array}{cc} 0 & l=0,\\ -\frac{\partial F}{\partial\kappa_1}(\bm{\kappa}_0)\frac{l(l+n+1)}{R^2} & l\in\mathbb{N},\end{array}\right.
\end{equation*}
with eigenfunctions
\begin{equation*}
u_{l,p}=\left\{\begin{array}{cc} Y_{0,1} & l=p=0,\\ Y_{l+1,p} & l\in\mathbb{N}\cup\{0\},1\leq p\leq M_{l+1}.\end{array}\right.
\end{equation*}
The multiplicity of the $0$ eigenvalue is then $M_1+1=n+2$.
$\Box$

In what follows, we set $P$ to be the projection from $C^{\alpha}\left(\SR\right)$ onto the $\lambda=0$ eigenspace given by
\begin{equation*}
Pu:=\sum_{p=0}^{n+1}\left\langle u,u_{0,p}\right\rangle u_{0,p},
\end{equation*}
where we use $\left\langle \cdot,\cdot\right\rangle$ to denote the $L^2$ inner product on $C^{2+\alpha}\left(\SR\right)$. Because $\partial G(0)$ is self adjoint with respect to this inner product, it is clear that $P\partial G(0)u=\partial G(0)Pu=0$ for every $u\in C^{2+\alpha}\left(\SR\right)$. Due to this $C^{2+\alpha}\left(\SR\right)$ can be split into the subspaces $X^c=P\left(C^{2+\alpha}\left(\SR\right)\right)$ and $X^s=(I-P)\left(C^{2+\alpha}\left(\SR\right)\right)$, called the \textit{center subspace} and \textit{stable subspace} respectively. We are now in a position to apply Theorem 9.2.2 in \cite{Lunardi95}.

\begin{theorem}
For any $l\in\mathbb{N}$ there is a function $\gamma\in C^{l-1}\left(X^c,X^s\right)$ such that $\gamma^{(l-1)}$ is Lipschitz continuous, $\gamma(0)=\partial\gamma(0)=0$ and $\mathcal{M}^c=\textrm{graph}(\gamma)$ is a locally invariant manifold for the equation (\ref{GraphLin}) of dimension $n+2$.
\end{theorem}

Note that by \textit{locally invariant} it is meant that there exists a neighbourhood of zero in $\Lambda\subset X^c$ such that if $\rho_0\in\textrm{graph}\left(\gamma|_{\Lambda}\right)$ then the solution to (\ref{GraphLin}) is in $\textrm{graph}\left(\gamma|_{\Lambda}\right)$ for all time or until $P\rho(t)\notin\Lambda$. We now set
\begin{equation*}
\mathcal{S}:=\left\{\rho\in C^{2+\alpha}\left(\SR\right):\textrm{graph}(\rho)\textrm{ is a sphere}\right\}.
\end{equation*}

\begin{lemma}
$\mathcal{M}^c$ coincides with the set $\mathcal{S}$ on a neighbourhood of zero.
\end{lemma}

\proof
Firstly note that due to Theorem 2.4 in \cite{Simonett95}, Theorem 2.3 in \cite{Iooss92} can be applied to conclude $\mathcal{M}^c$ contains all equilibria of (\ref{GraphLin}) with $P\rho_0\in\Lambda$. Also note, the center manifold, $\mathcal{M}^c$, is defined differently in \cite{Iooss92} as compared to \cite{Lunardi95}, however they can be seen to be equal on $\Lambda$. The rest of the proof follows from \cite{Escher98A}, where $\rho\in\mathcal{S}$ is formulated in terms of the eigenfunctions:
\begin{equation*}
\rho(\bm{z})=\sum_{p=1}^{n+1}z_pu_{0,p}-R+\sqrt{\left(\sum_{p=1}^{n+1}z_pu_{0,p}\right)^2+(R+z_0)^2-\sum_{p=1}^{n+1}z_p^2},
\end{equation*}
with $(z_1,\ldots,z_{n+1})\in\Rnone$ being the centre of the sphere and $z_0:=R'-R$ the difference between its radius and the radius of $\SR$. This map is smooth on a neighbourhood $U$ of $0\in\mathbb{R}^{n+2}$ and its derivative at zero is given by
\begin{equation*}
\partial\rho(0)\bm{z}=\sum_{p=0}^{n+1}z_pu_{0,p},\ \ \ \bm{z}\in\mathbb{R}^{n+2}.
\end{equation*}
The map taking $\bm{z}$ to the coordinates of $P\rho(\bm{z})$ with respect to the basis $u_{0,p}$, $0\leq p\leq n+1$, is then found to have derivative at zero equal to the identity and hence is a diffeomorphism from $U$ onto its image, possibly making $U$ smaller. This means that the projection of $\mathcal{S}|_{U}:=\{\rho(\bm{z}):\bm{z}\in U\}$ is an open neighbourhood of $0\in X^c$. This can be made to coincide with $\Lambda$ (after possible renaming) and since $\mathcal{S}|_{U}\subset\mathcal{M}^c$ (by the first remark of this proof) we conclude that $\mathcal{S}$ and $\mathcal{M}^c$ coincide locally.
$\Box$

We now prove the main result.
\proofoftheoremone
By Proposition 9.2.4 in \cite{Lunardi95} for every $\omega\in(0,-\lambda_1)$ there is a constant $C(\omega)>0$ and a neighbourhood, $O_4$, of $0\in C^{2+\alpha}\left(\SR\right)$ such that if $\rho_0\in O_4$ there exists $\bar{x}\in \Lambda$ such that
\begin{equation*}
\left\|P\rho(t)-\bar{x}\right\|_{C^{2+\alpha}\left(\SR\right)}+\left\|(I-P)\rho(t)-\gamma(\bar{x})\right\|_{C^{2+\alpha}\left(\SR\right)}\leq C(\omega)e^{-\omega t}\left\|(I-P)\rho_0-\gamma(P\rho_0)\right\|_{h^{C+\alpha}\left(\SR\right)},\ \ \ \forall t\geq0.
\end{equation*}
Here we have used that for every $x\in\Lambda$ the function $x+\gamma(x)$ defines a sphere and hence is a stationary solution to (\ref{GraphFlow}). This proves that $\rho(t)$ converges to an element of $\mathcal{M}^c|_{\Lambda}$, which by the above lemma is a sphere.
$\Box$

\begin{corollary}
Let $\Omega_0$ be a graph over a sphere with height $\rho_0$ such that the solution, $\rho(t)$, to the flow (\ref{GraphFlow}) with initial condition $\rho_0$ exists for all time and converges to zero. Suppose further that $\left.\frac{\partial F}{\partial \kappa_i}\right|_{\kappa_{\rho(t)}}>0$ for all $t\in[0,\infty)$ and $i=1,\ldots,n$. Then there exists a neighbourhood, $O$, of $\rho_0$ in $C^{2+\alpha}\left(\SR\right)$, $0<\alpha<1$, such that for every $u_0\in O$ the solution to (\ref{GraphFlow}) with initial condition $u_0$ exists for all time and converges to a function near zero whose graph is a sphere.
\end{corollary}

\proof
This follows by the same arguments given in \cite{Guenther02}. First we set $U\subset C^{2+\alpha}\left(\SR\right)$ to be the neighbourhood of zero given in Theorem \ref{MainR}. Since $\rho(t)$ converges to zero in the $C^{2+\alpha}$-topology there exists a time $T$ such that $\rho(T)\in U$ and as $U$ is open there exists an open ball $B_{\epsilon}(\rho(T))\subset U$ of radius $\epsilon$ centred at $\rho(T)$. The condition that $\left.\frac{\partial F}{\partial \kappa_i}\right|_{\kappa_{\rho(t)}}>0$ for all $t\in[0,\infty)$ and $i=1,\ldots,n$ ensures that the operator $L_{\rho}F(\bm{\kappa}_{\rho})$ is elliptic around the point $\rho(t)$ for every $t\in[0,\infty)$ (see \cite{Andrews94}). By taking care of the linearisation of the global term in the same way as in the proof of Theorem \ref{MainR} the linear operator $\partial G(\rho(t))$ can be seen to be sectorial of all $t\in[0,T]$, and hence in a neighbourhood of each point. By Theorem 8.4.4 in \cite{Lunardi95} the flow depends continuously on the initial condition in a neighbourhood of $\rho_0$. Therefore there exists a ball $B_{\delta}(\rho_0)$ such that if $u_0\in B_{\delta}(\rho_0)$ then the solution, $u(t)$, to (\ref{GraphFlow}) with initial condition $u_0$ exists for $t\in[0,T]$ and $u(T)\in B_{\epsilon}(\rho(T))$. Since $u(T)$ is in $U$, by Theorem \ref{MainR} the solution to (\ref{GraphFlow}) with initial condition $u(T)$ converges to a function near zero that defines a sphere. By uniqueness of the flow we get the result.
$\Box$

\bibliographystyle{plain}
\bibliography{CompleteWorkNew}

\end{document}